\documentclass[10pt,leqno,twoside,draft]{article}
\usepackage{amsmath}
\usepackage{amsfonts}
\usepackage{amsthm}
\usepackage{indentfirst}

\theoremstyle{plain}
\newtheorem{theorem}{\indent\rm T\,h\,e\,o\,r\,e\,m\;}[section]

\newtheorem{proposition}{\indent\rm P\,r\,o\,p\,o\,s\,i\,t\,i\,o\,n\;}[section]

\theoremstyle{definition}
\newtheorem{definition}{\indent\rm D\,e\,f\,i\,n\,i\,t\,i\,o\,n\;}[section]

\theoremstyle{remark}
\newtheorem{remark}{\indent\rm R\,e\,m\,a\,r\,k\;}[section]

\makeatletter                                                  %
\renewcommand*{\@seccntformat}[1]{
  \csname the#1\endcsname\;-                                   %
}                                                              %
\renewcommand{\section}{\@startsection{section}{1}{0mm}        %
   {1.5\baselineskip}
   {1\baselineskip}
   {\indent\normalfont\normalsize\bfseries}
   }                                                           %
\renewcommand*{\@seccntformat}[1]{
  \normalfont\bfseries\csname the#1\endcsname\;-               %
}                                                              %
\renewcommand\subsection{\@startsection                        %
  {subsection}{2}{0mm}
  {1.5\baselineskip}
  {1\baselineskip}
  {\indent\normalfont\normalsize\itshape}}
\renewcommand*{\@seccntformat}[1]{
  \normalfont\bfseries\csname the#1\endcsname\;-               %
}                                                              %
\renewcommand\subsubsection{\@startsection                     %
  {subsubsection}{2}{0mm}
  {1.5\baselineskip}
  {1\baselineskip}
  {\indent\normalfont\normalsize\texttt}}
\makeatother                                                   %
\begin{document}
\thispagestyle{empty}

\begin{center}
{\sc\large Elena Angelini} 
\end{center}
\vspace {1.5cm}

\centerline{\large{\textbf{On complex and real identifiability of tensors}}}

\renewcommand{\thefootnote}{\fnsymbol{footnote}}

\footnotetext{
This research was partially supported by the Italian GNSAGA-INDAM and by the Italian PRIN2015 - Geometry of Algebraic Varieties (B16J15002000005)}

\renewcommand{\thefootnote}{\arabic{footnote}}
\setcounter{footnote}{0}

\vspace{1,5cm}
\begin{center}
\begin{minipage}[t]{10cm}

\small{ \noindent \textbf{Abstract.} We report about the state of the art on complex and real generic identifiability of tensors, we describe some of our recent results obtained in \cite{ABC} and we present perspectives on the subject.  
\medskip

\medskip

\noindent \textbf{Mathematics~Subject~Classification~(2010):}
14N05, 15A69, 14P05.

}
\end{minipage}
\end{center}

\bigskip

\section{Introduction}

Identifiability is a very special property for tensors. This interesting topic has been extensively investigated starting from the $XIX^{th}$ century, related to the Waring problem \cite{Hi,K,R,S,Sy,We}, however the theory is far from being complete. The notion of identifiability has many important applications, that go beyond Algebraic Geometry: for example, in engineering, with the Blind Source Separation problem \cite{AGHKT}, in chemistry, when one deals with mixtures of fluorophores \cite{AD}, or in statistics \cite{AMR}. 

The paper is organized as follows. In section $ 2 $ we introduce preliminary material. Section $ 3 $ is devoted to the complex case, with particular emphasis on symmetric, partially-symmetric and skew-symmetric tensors, which are three significant special cases. In section $ 4 $ we focus on the real case, describing some of our recent results obtained in \cite{ABC}. Finally in section $ 5 $ we present open problems and perspectives on the subject. 

The main purpose of this note is to report about the state of the art on generic identifiability of tensors, stressing, in a clear and simple way, analogies and differences between the complex and real case.

\section{Preliminaries}

In this section we recall basic definitions about tensor decomposition, rank and identifiability, mainly referring to the introductory handbook \cite{Lan}. Moreover we describe some properties of elliptic normal curves, that play a special role in the discovery of unidentifiable cases. 

\subsection{Tensor decomposition and rank}

Let $\mathbb{F}$ be either the complex field $ \mathbb{C} $ or the real field $ \mathbb{R} $, let $ n_{1}, \ldots, n_{d} \in \mathbb{N} $ such that $ n_{1} \leq \ldots \leq n_{d} $ and let $ \mathbb{F}^{n_{1}} \otimes \ldots \otimes \mathbb{F}^{n_{d}} $ be the space of tensors of type $ n_{1}, \ldots, n_{d} $ over $\mathbb{F}$. Let $ T \in\mathbb{F}^{n_{1}} \otimes \ldots \otimes \mathbb{F}^{n_{d}} $.
\begin{definition}
A \emph{tensor decomposition} of $ T $ is given by vectors $ {\bold{v}}_{i}^{j} \in\mathbb{F}^{n_{j}} $, with $ i \in\{1, \ldots, k\} $ and $ j \in \{1, \ldots, d\} $, such that
\begin{equation}\label{eq:tensordec}
T = \displaystyle{\sum_{i=1}^{k} {\bold{v}}_{i}^{1} \otimes \ldots \otimes {\bold{v}}_{i}^{d}}. 
\end{equation}
\end{definition}
\begin{definition}
The minimal $ k $ appearing in (\ref{eq:tensordec}) is the \emph{rank} of $ T $ \emph{over} $ \mathbb{F} $. 
\end{definition}

\begin{remark}
Any $ {\bold{v}}_{i}^{1} \otimes \ldots \otimes {\bold{v}}_{i}^{d} $ has rank $ 1 $ over $ \mathbb{F} $.
\end{remark}

\begin{definition}
A \emph{typical rank} over $ \mathbb{F} $ for tensors of type $ n_{1}, \ldots, n_{d} $ is any $ k $ such that the set of tensors having rank $ k $ has positive Euclidean measure.
\end{definition}

According to section 5.2.1 of \cite{Lan}, we introduce the following:

\begin{definition}
There exists a unique typical rank for tensors of type $ n_{1},\dots, n_{d} $  over $ \mathbb{C} $, which we call the \emph{generic rank} for that space of tensors. 
\end{definition}

\begin{remark}
The \emph{expected generic rank} for $ \mathbb{C}^{n_{1}} \otimes \ldots \otimes \mathbb{C}^{n_{d}} $ is 
$$ k_{g} = \left\lceil \frac{\prod_{i=1}^{d}n_{i}}{1+\sum_{i=1}^{d}(n_{i}-1)} \right\rceil. $$
\end{remark}

According to \cite{COV}, we give the following:
\begin{definition}
Tensors of type $ n_{1}, \dots , n_{d} $ over $ \mathbb{C} $ and rank $ k < k_{g} $ are of \emph{sub-generic rank}.
\end{definition}

\begin{remark}\label{rem:c&r}
It is well known that it is possible to have more than one typical rank for tensors of type $ n_{1}, \ldots, n_{d} $ over $ \mathbb{R} $. The smallest typical rank over $ \mathbb{R} $ coincides with the generic rank over $ \mathbb{C} $ \cite{BT}. Any rank over $ \mathbb{R} $ between the minimal typical rank and the maximal typical rank is also typical and it is an open problem to determine an expression for the expected maximal typical rank over $ \mathbb{R} $ in the general case. Partial results have been obtained for symmetric tensors  \cite{BBO}.
\end{remark}

\subsection{Identifiability}\label{subsec:id}

\begin{definition}
A rank-$ k $ tensor $ T \in \mathbb{F}^{n_{1}} \otimes \ldots \otimes \mathbb{F}^{n_{d}} $ is \emph{identifiable over} $ \mathbb{F} $ if the presentation (\ref{eq:tensordec}) is unique up to a permutation of the summands and scaling of the vectors. 
\end{definition} 
The above definition extends as follows:
\begin{definition}
The set of tensors of type $ n_{1}, \ldots, n_{d} $ and rank $ k $ over $\mathbb{F}$  is \emph{generically identifiable over} $ \mathbb{F} $  if identifiability over $ \mathbb{F} $ holds on a Zariski dense open subset of the variety of tensors of rank $ \leq k $. 
\end{definition}

\subsection{Elliptic normal curves}\label{subsec:enc}

\indent Let $ \mathbb{P}^{n} = {\mathbb{P}}_{\mathbb{C}}^{n} $ be the $ n $-dimensional complex projective space.
\begin{definition}
An \emph{elliptic normal curve} $ \mathcal{C} \subset \mathbb{P}^{n} $ is a smooth curve of genus $ 1 $ and degree $ n+1 $ that is contained in no hyperplane.
\end{definition} 
For any such curve $ \mathcal{C} \subset \mathbb{P}^{2k-1} $ of even degree $ n+1 = 2k $, it is known that $ k-1 $ is the minimal dimension of a $ \mathbb{P}^{h-1} $ $ h $-secant to $ \mathcal{C} $ which contains the general point $ P \in \mathbb{P}^{2k-1} $. In particular the following holds:

\begin{proposition}[Chiantini-Ciliberto, 2006, \cite{CC}]\label{prop:CC}
Under the above assumptions, the number of $ k $-secant $ \mathbb{P}^{k-1} $ to $ \mathcal{C} $ passing through $ P $ is $ 2 $.
\end{proposition}

\begin{remark}
As a consequence of Proposition \ref{prop:CC}, we can prove the existence of two decompositions over $ \mathbb{C} $ in some cases. Precisely, up to now, the only known examples are those listed in Table $ 4 $ of section $ 4 $. Indeed, let $ \mathbb{T} $ be the space of tensors object of study and $  \mathbb{T}_{1} \subset \mathbb{P}(\mathbb{T}) $ the corresponding variety of rank-$1$ tensors over $ \mathbb{C} $. We choose a general tensor $ T \in \mathbb{P}(\mathbb{T}) $ of rank $ k $ over $ \mathbb{C} $ and one of its complex decompositions, i.e.\ $T_{1}, \ldots, T_{k} \in \mathbb{T}_{1} $ such that $ T = \sum_{i=1}^{k} T_{i} $. We can show that there exists a unique irreducible elliptic normal curve $ \mathcal{C} $ of even degree $ 2k $ entirely contained in $ \mathbb{T}_{1} $, passing through the summands $ T_{1}, \ldots, T_{k} $ and such that $  T $ belongs to the odd dimensional projective space $ \mathbb{P}^{2k-1} $ spanned by $ \mathcal{C} $. By Proposition \ref{prop:CC}, $ T $ belongs to two $\mathbb{P}^{k-1}$ $k$-secant to $ \mathcal{C} $, which provide the two complex decompositions of $ T $. We notice that, for the cases of sub-generic rank, the elliptic normal curve is the \emph{$ k $-contact locus} of $  \mathbb{T}_{1} $ (see \cite{CC} for more details on this concept). As will be explained in our forthcoming paper \cite{ABC1}, for the last example of Table $ 4 $ it is not possible to define a contact locus, since in this case the $k$-secant variety of $ \mathbb{T}_{1} $ fills the ambient space, so that the curve represents a ``known" subvariety of $ \mathbb{T}_{1} $ bounding the decompositions of $ T $.
\end{remark}

\section{Complex identifiability}

Generic identifiability over $ \mathbb{C} $ has been largely investigated, in particular we refer to \cite{BCV,BC,BCO,CO,COV,COV1,K,S}. 

The first interesting result concerns the case $ d = 3 $:
\begin{theorem}[Kruskal, 1977, \cite{K}]\label{thm:K}
The set of tensors of type $ n_{1},n_{2},n_{3} $ and rank $ k $ over $ \mathbb{C} $ is generically identifiable over $ \mathbb{C} $ if
$$ k \leq \frac{1}{2} (\min(n_{1},k) + \min(n_{2},k) + \min(n_{3},k) -2). $$ 
\end{theorem}
Theorem \ref{thm:K} has been improved for $ n_{1} = n_{2} = n_{3} = n $ odd as follows:
\begin{theorem}[Strassen, 1983, \cite{S}]
The set of tensors of type $ n, n, n $ with $ n $ odd and rank $ k $ over $ \mathbb{C} $ is generically identifiable over $ \mathbb{C} $ if
$$ k \leq \left\lfloor\frac{n^3}{3n-2} \right\rfloor - n. $$ 
\end{theorem}
Furthermore we have the following:
\begin{theorem}[Bocci-Chiantini-Ottaviani, 2014, \cite{BCO}]
The set of tensors of type $ n_{1},n_{2},n_{3} $ with $ 3 < n_{1} \leq n_{2} \leq n_{3} $ and rank $ k $ over $ \mathbb{C} $ is generically identifiable over $ \mathbb{C} $ if
$$ k \leq \frac{n_{1}n_{2}n_{3}}{n_{1}+n_{2}+n_{3}-2}-n_{3}. $$ 
\end{theorem}

Admitting a unique minimal decomposition is a quite rare phenomenon for tensors of generic rank over $ \mathbb{C} $. As an example, in the \emph{symmetric} setting, i.e. \ when $ n_{1} = \ldots = n_{d} = n $ and we deal with Sym$^{d}\mathbb{C}^n  \subset (\mathbb{C}^n)^{\otimes d}$, the complex identifiable cases of tensors of generic rank are only the one classically known described in Table $ 1 $, as stated in \cite{GM}. We would like to emphasize that the first and third examples are due to Sylvester, while the second one to Hilbert.
\begin{table}[h]
\begin{center}
\begin{tabular} {c |c |c |c }
{\bf Space of tensors} & {\bf Feature} & {\bf Rank}  & {\bf Ref.} \\
$(\mathbb{C}^2)^{\otimes 2t+1}$ & symmetric & $t+1$ & \cite{Sy} \\
$(\mathbb{C}^3)^{\otimes 5}$ & symmetric & $7$ & \cite{Hi} \\
$(\mathbb{C}^4)^{\otimes 3}$ & symmetric & $5$ & \cite{Sy} \\
\end{tabular}\caption{Symmetric generically identifiable cases of generic rank over $ \mathbb{C} $}
\end{center}
\end{table}

\noindent In the \emph{partially-symmetric} case, that is when we work with vectors of symmetric tensors in the same variables, searching for simultaneous decompositions, the complete list of identifiable cases of generic rank is not known. The discovered ones, up to now, are collected in Table $ 2 $. Besides of Veronese's result, which goes back to 1880, the second and fourth examples are classical too, respectively, due to Weierstrass and Roberts. For a more detailed discussion on simultaneous decompositions of vectors of symmetric tensors and identifiability, we refer to our papers \cite{A} and \cite{AGMO}.  
\begin{table}[h]
\begin{center}
\begin{tabular} {c |c |c |c }
{\bf Space of tensors} & {\bf Feature} & {\bf Rank}  & {\bf Ref.} \\
$\bigoplus_{j=1}^{r} (\mathbb{C}^2)^{\otimes d_{j}}, \, d_{1}+1 \geq k $ & part.-symm. & $ \left \lceil {\frac{1}{1+r}}\sum_{j=1}^{r}{1+d_{j} \choose d_{j}} \right \rceil $ & \cite{CR} \\
$((\mathbb{C}^n)^{\otimes 2})^{\oplus 2} $ & part.-symm. & $n+1$ & \cite{We} \\
$((\mathbb{C}^3)^{\otimes 2})^{\oplus 4}$ & part.-symm. & $4$ & Veronese \\
$(\mathbb{C}^3)^{\otimes 2} \oplus (\mathbb{C}^3)^{\otimes 3} $ & part.-symm. & $4$ & \cite{R} \\
$((\mathbb{C}^3)^{\otimes 3})^{\oplus 2} \oplus (\mathbb{C}^3)^{\otimes 4} $ & part.-symm. & $7$ & \cite{AGMO} \\
\end{tabular}\caption{Part.-symm. generically identifiable cases of generic rank over $ \mathbb{C} $}
\end{center}
\end{table}

On the other side, generic identifiability over $ \mathbb{C} $ is expected for tensors of sub-generic rank over $ \mathbb{C} $. In this direction, by means of an algorithm based on the so-called \emph{Hessian criterion}, in \cite{COV} it is proved generic identifiability over $ \mathbb{C} $ for a large number of spaces of tensors in sub-generic rank cases. For the symmetric case (with $ d \geq 3 $), in \cite{COV1} it is proved that there are no exceptions besides the ones appearing in Table $ 3 $.
\begin{table}[h]
\begin{center}
\begin{tabular} {c |c |c |c }
{\bf Space of tensors} & {\bf Feature} & {\bf Rank}  & {\bf Ref.} \\
$(\mathbb{C}^3)^{\otimes 6}$ & symmetric & $9$  & \cite{AC,CC} \\
$(\mathbb{C}^4)^{\otimes 4}$ & symmetric & $8$  & \cite{CC,M}\\
$(\mathbb{C}^6)^{\otimes 3}$ & symmetric & $9$  & \cite{COV1,RV}\\
\end{tabular}\caption{Symmetric generically unidentifiable cases of sub-generic rank over $ \mathbb{C} $}
\end{center}
\end{table}

\noindent An account of this topic in the partially-symmetric setting can be found in our forthcoming paper \cite{A}.

Concerning generic identifiability over $ \mathbb{C} $ of \emph{skew-symmetric} cases,  i.e.\ when $ n_{1} = \ldots = n_{d} = n $ and it is investigated $ \bigwedge^d \mathbb{C}^n \subset (\mathbb{C}^n)^{\otimes d}$, we refer to \cite{BV} (Theorem 1.1).

Besides the cases described above, generic identifiability over $ \mathbb{C} $ does not hold anytime the projective algebraic variety $ \mathbb{T}_{1} $ of rank-$1$ tensors object of study (Segre, Veronese and Grassmann variety, respectively, for the general, the symmetric and the skew-symmetric case) is \emph{$ k $-defective} (we refer to \cite{Lan} for the basics about secant varieties and the defectivity problem and to \cite{AOP,AH,BDdG} for an account on Segre, Veronese and Grassmann defective varieties), in which case we deal with rank-$k$ tensors admitting infinitely many complex decompositions. 

\section{Real identifiability}\label{sec:recres}

Recent interest has been devoted to the real case \cite{COV2,CLQ,DDL}, very useful in applications. 

If a set of tensors, of fixed type and rank over $\mathbb{C}$,  is generically identifiable over $ \mathbb{C} $, then, necessarily, the same set of tensors seen over $\mathbb{R}$ is generically identifiable over $ \mathbb{R} $. Equivalently, with the above assumption, the unique complex decomposition of the general real tensor of this type and rank over $\mathbb{R}$ is completely real.

At this point a natural question arises: if identifiability over $ \mathbb{C} $ fails, what happens to identifiability over $\mathbb{R}$? For example, one may wonder if there exist real tensors admitting several decompositions in terms of complex rank-$ 1 $ summands but only one of them is a decomposition over $ \mathbb{R} $.

In \cite{ABC} we answered this question in some cases. 

The main result is for tensors having two complex decompositions explained through elliptic normal curves:

\begin{theorem}[Angelini-Bocci-Chiantini, 2017, \cite{ABC}]\label{thm:ABC}
Let $ \mathbb{T} $ be any space of tensors in Table $ 4 $ of fixed feature and $ k $ the rank over $ \mathbb{C} $ under investigation. 
\begin{table}[h]
\begin{tabular} {c |c |c |c }
{\bf Space of tensors} & {\bf Feature} & {\bf Rank}  & {\bf Ref.} \\
$(\mathbb{C}^2)^{\otimes 5}$ &  & $5$ (sg) & \cite{BC} \\
$(\mathbb{C}^3)^{\otimes 6}$ & symmetric & $9$ (sg) & \cite{AC,CC} \\
$(\mathbb{C}^4)^{\otimes 3}$ &  & $6$ (sg) & \cite{CO} \\
$(\mathbb{C}^4)^{\otimes 4}$ & symmetric & $8$ (sg) & \cite{CC,M}\\
$(\mathbb{C}^6)^{\otimes 3}$ & symmetric & $9$ (sg) & \cite{COV1,RV}\\
$(\mathbb{C}^{10})^{\otimes 3}$ & skew-symmetric & $5$ (sg) & \cite{BV} \\ 
$((\mathbb{C}^{3})^{\otimes 3})^{\oplus 3} \oplus ((\mathbb{C}^{3})^{\otimes 2})^{\oplus m}, \, m \geq 0 $ & partially-symmetric & $6$ (g) & \cite{ABC1,L} \\ 
\end{tabular}\caption{Gen. unidentifiable cases of sub-generic and generic rank over $ \mathbb{C} $}
\end{table} 

\noindent Then there exist non-trivial Euclidean open subsets $ U_{1}, U_{2}, U_{3} $ of the variety of real tensors in $ \mathbb{T} $ of rank $ \leq k $ over $ \mathbb{C} $, whose elements have two complex decompositions and:
\begin{itemize}
\item[$\bullet$] $ \forall \, T \in U_{1} $, only one decomposition is real; 
\item[$\bullet$] $ \forall \, T \in U_{2} $, both decompositions are real;
\item[$\bullet$] $ \forall \, T \in U_{3} $, both decompositions are not real.
\end{itemize}

\end{theorem}

In column $ 2 $ of Table $ 4 $, when not declared, we intend any tensor of the corresponding space. The last row is devoted to the partially-symmetric case of vectors of symmetric tensors with three ternary cubics and an arbitrary number of conics. This example will be extensively treated in our forthcoming paper \cite{ABC1}.

The proof of Theorem \ref{thm:ABC}, the details of which we refer to \cite{ABC}, is based on a study of real elliptic normal curves of even degree $ 2k $ in projective spaces of odd dimension $ 2k-1 $, through techniques of projective geometry. In particular we first describe the case of quartics showing the following:
 
\begin{proposition}[Angelini-Bocci-Chiantini, 2017, \cite{ABC}]\label{prop:quartics}
Let $ \mathcal{C} \subset \mathbb{P}^{3} $ be an irreducible real elliptic normal quartic. Then there exist $ B_{1}, B_{2}, B_{3}, B_{4} \subset {\mathbb{P}}_{\mathbb{R}}^{3} $ non-trivial open balls entirely composed of points $ P $ such that the two secant lines to $ \mathcal{C} $ through $ P $ intersect $ \mathcal{C} $, respectively, in $ 4 $ real points, $ 2 $ real and $ 2 $ not real conjugate points, $ 4 $ not real points pairwise conjugate (for $ B_{3} $ the conjugate points being on the same line and for $ B_{4} $ otherwise).
\end{proposition}

By induction on $ k $ we extend Proposition \ref{prop:quartics} to curves in higher dimensional spaces getting the following:

\begin{proposition}[Angelini-Bocci-Chiantini, 2017, \cite{ABC}]\label{prop:curves}
Let $ \mathcal{C} \subset \mathbb{P}^{2k-1} $ be an irreducible real elliptic normal curve of degree $ 2k $. Then, for any $ (a_{1}, a_{2})\in \mathbb{N}^2 $ such that $ 0 \leq 2a_{1}+2a_{2} \leq k $, there exists $ B_{a_{1},a_{2}}\subset {\mathbb{P}}_{\mathbb{R}}^{2k-1} $ non-trivial open ball entirely composed of points $ P $ such that each of the two secant spaces $ \Pi_{i}(P) \cong \mathbb{P}^{k-1} $, $i \in \{1,2\}$, to $ \mathcal{C} $ through $ P $ intersects $ \mathcal{C} $ in k-2$a_{i}$ real points and 2$a_{i}$ not real points.
\end{proposition}

We notice that, with the notation introduced in Proposition \ref{prop:curves}, if one of the $ a_{i} $'s equals $ 0 $ and the other one is different from $ 0 $, then we get real points $ P $ admitting a unique real decomposition with respect to $ \mathcal{C} $. If $ a_{1} = a_{2} = 0 $, then the two decompositions are real. Finally, if the $ a_{i} $'s are both different from $ 0 $, then the two decompositions are complex but not real. In particular, the first case will lead to real identifiability.

Proposition \ref{prop:curves} can be generalized to families of irreducible real elliptic normal curves. This step allows us to pass from curves to tensors, obtaining the proof of Theorem \ref{thm:ABC}.

\section{Open problems and perspectives}\label{sec:persoppbs}

We conclude this note by presenting a short list of open problems that have arisen mainly from discussions within the seminars I held on the subject.\\
\indent $1)$ Classifying spaces of tensors, of fixed type and rank, such that the general tensor has two complex decompositions, due to the presence of elliptic normal curves, is, at the moment, an open problem. 

$2)$ According to Remark \ref{rem:c&r}, from Theorem \ref{thm:ABC} we immediately get some information about the rank over $ \mathbb{R} $ and identifiability over $ \mathbb{R} $ of tensors in $ U_{1}, U_{2}, U_{3} $. Indeed, in $ U_{1} $ the rank over $ \mathbb{R} $ equals the rank over $ \mathbb{C} $ and real identifiability holds. In $ U_{2} $ the same condition on the ranks holds but identifiability over $ \mathbb{R} $ fails. In $ U_{3} $ the rank over $ \mathbb{R} $ is strictly greater than the rank over $ \mathbb{C} $ and we can't say anything about real identifiability. Therefore we can ``recover" identifiability only in $ U_{1} $. At this point it is reasonable to ask how large the three open sets of Theorem \ref{thm:ABC} are and, till now, also determining an estimate of their measure represents an interesting open problem. 

$3)$ It is not known completely how identifiability over $ \mathbb{R} $ behaves when identifiability over $ \mathbb{C} $ fails for reasons other than the presence of elliptic normal curves. In this direction, in \cite{ABC} we showed that, if the general tensor of rank $ k $ over $ \mathbb{C} $ admits infinitely many decompositions, then there are no Euclidean open sets of tensors in which real identifiability holds. 

$4)$ Several situations seem to occur in the cases with a finite number, greater than $ 1 $, of complex decompositions. For instance, let us consider tensors with complex sub-generic rank $ k = n_{d}-1 $ in \emph{almost unbalanced} spaces $ \mathbb{C}^{n_{1}} \otimes \ldots \otimes \mathbb{C}^{n_{d}} $ (i.e.\ $ n_{1} \leq \ldots \leq n_{d} $ and $ n_{d}-1 = \prod_{i=1}^{d-1}(n_{i}) - \sum_{i=1}^{d-1}(n_{i}-1) $: it is known, by \cite{BCO}, that they have $ {D}\choose{k} $ complex decompositions, with $ D = \deg(\mathbb{P}^{n_{1}-1} \times \ldots \times \mathbb{P}^{n_{d-1}-1} ) $, and we proved in \cite{ABC}  that, whenever $ D - n_{d}+1 $ is odd, then there are no Euclidean open sets where identifiability over $ \mathbb{R} $ holds, while, when $ D - n_{d}+1 $ is even, we believe in the existence of such open sets. In this sense, this is true when $ d = n_{1} = 3, n_{2} = 5, n_{3} = 10 $, see Example 5.2 of \cite{ABC}. Concerning the cases of complex generic rank, we can apply the computer-aided procedure introduced in \cite{AGMO} and based on \emph{homotopy continuation techniques} and \emph{monodromy loops}, to produce examples where real identifiability is recovered at least in Euclidean open sets (see Table $ 5 $ for the cases analyzed with this method in our papers \cite{ABC} and \cite{A}). 

\begin{table}[h]
\begin{center}
\begin{tabular} {c |c |c |c |c }
{\bf Space of tensors} & {\bf Feature} & {\bf Rank}  & {\bf Dec.} & {\bf Ref.} \\
$(\mathbb{C}^3)^{\otimes 7}$ & symmetric & $12$ (g) & 5 & \cite{RS} \\
$(\mathbb{C}^3)^{\otimes 8}$ & symmetric & $15$ (g) & 16 & \cite{RS} \\
$((\mathbb{C}^{3})^{\otimes 3})^{\oplus 3} \oplus (\mathbb{C}^{3})^{\otimes 2} $ & partially-symmetric & $6$ (g) & 2 & \cite{ABC1} \\
\end{tabular}\caption{Generically unidentifiable cases analyzed through computations}
\end{center}
\end{table}

\vspace{0.5cm} \indent {\it
A\,c\,k\,n\,o\,w\,l\,e\,d\,g\,m\,e\,n\,t\,s.\;} This note, based on joint works with C. Bocci and L. Chiantini, arises partially from the conference ``Real identifiability and complex identifiability" given by the author within the cycle of study seminars on ``Algebraic Geometry and Tensors", held among Bologna, Ferrara, Firenze and Siena in 2016-2017. The author would like to thank all the participants for fruitful and stimulating discussions.

\bigskip
\begin{center}

\end{center}

\bigskip
\bigskip
\begin{minipage}[t]{10cm}
\begin{flushleft}
\small{
\textsc{Elena Angelini}
\\*University of Siena,
\\*Dipartimento di Ingegneria dell'Informazione e Scienze Matematiche
\\*Via Roma 56
\\* Siena, 53100, Italy
\\*e-mail: elena.angelini@unisi.it
}
\end{flushleft}
\end{minipage}

\end{document}